\documentclass[a4paper,10pt]{amsart}

\usepackage{amsfonts}

\usepackage{graphicx}
\usepackage{amsmath}
\usepackage{amssymb}
\usepackage[all]{xypic}
\usepackage[all]{xy}

\begin{document}
\input xy
\xyoption{all}

\renewcommand{\mod}{\operatorname{mod}\nolimits}
\newcommand{\proj}{\operatorname{proj}\nolimits}
\newcommand{\rad}{\operatorname{rad}\nolimits}
\newcommand{\soc}{\operatorname{soc}\nolimits}
\newcommand{\ind}{\operatorname{ind}\nolimits}
\newcommand{\Top}{\operatorname{top}\nolimits}
\newcommand{\ann}{\operatorname{Ann}\nolimits}
\newcommand{\id}{\operatorname{id}\nolimits}
\newcommand{\Mod}{\operatorname{Mod}\nolimits}
\newcommand{\End}{\operatorname{End}\nolimits}
\newcommand{\Ob}{\operatorname{Ob}\nolimits}
\newcommand{\Ht}{\operatorname{Ht}\nolimits}
\newcommand{\cone}{\operatorname{cone}\nolimits}
\newcommand{\rep}{\operatorname{rep}\nolimits}
\newcommand{\Ext}{\operatorname{Ext}\nolimits}
\newcommand{\Hom}{\operatorname{Hom}\nolimits}
\renewcommand{\Im}{\operatorname{Im}\nolimits}
\newcommand{\Ker}{\operatorname{Ker}\nolimits}
\newcommand{\Coker}{\operatorname{Coker}\nolimits}
\renewcommand{\dim}{\operatorname{dim}\nolimits}
\newcommand{\Ab}{{\operatorname{Ab}\nolimits}}
\newcommand{\Coim}{{\operatorname{Coim}\nolimits}}
\newcommand{\pd}{\operatorname{pd}\nolimits}
\newcommand{\sdim}{\operatorname{sdim}\nolimits}
\newcommand{\add}{\operatorname{add}\nolimits}
\newcommand{\cc}{{\mathcal C}}
\newcommand{\ct}{{\mathcal T}}
\newcommand{\cd}{{\mathcal D}}
\newcommand{\ch}{{\mathcal H}}
\newcommand{\cu}{{\mathcal U}}
\newcommand{\cv}{{\mathcal V}}

\newtheorem{theorem}{Theorem}[section]
\newtheorem{acknowledgement}[theorem]{Acknowledgement}
\newtheorem{algorithm}[theorem]{Algorithm}
\newtheorem{axiom}[theorem]{Axiom}
\newtheorem{case}[theorem]{Case}
\newtheorem{claim}[theorem]{Claim}
\newtheorem{conclusion}[theorem]{Conclusion}
\newtheorem{condition}[theorem]{Condition}
\newtheorem{conjecture}[theorem]{Conjecture}
\newtheorem{corollary}[theorem]{Corollary}
\newtheorem{criterion}[theorem]{Criterion}
\newtheorem{definition}[theorem]{Definition}
\newtheorem{example}[theorem]{Example}
\newtheorem{exercise}[theorem]{Exercise}
\newtheorem{lemma}[theorem]{Lemma}
\newtheorem{notation}[theorem]{Notation}
\newtheorem{problem}[theorem]{Problem}
\newtheorem{proposition}[theorem]{Proposition}
\newtheorem{remark}[theorem]{Remark}
\newtheorem{solution}[theorem]{Solution}
\newtheorem{summary}[theorem]{Summary}
\newtheorem*{thma}{Theorem}

\title
{Lifting to cluster-tilting objects in higher cluster categories}

\author{Pin Liu }
\address{
 Department of Mathematics\\
  Sichuan University\\
  610064 Chengdu \\
  P.R.China}
\address{
 Department of Mathematics\\
  Southwest Jiaotong University\\
  610031 Chengdu \\
  P.R.China}
\email{
\begin{minipage}[t]{5cm}
pinliu@yahoo.cn \\
\end{minipage}
}


\begin{abstract}  Let $d>1$  be a positive integer. In this note, we consider the $d$-cluster-tilted algebras, the endomorphism algebras
of $d$-cluster-tilting objects in $d$-cluster categories. We show
that a tilting module over such an algebra  lifts to a
$d$-cluster-tilting object in this $d$-cluster category.
\end{abstract}

\maketitle

\section{Introduction}
Let $k$ be an algebraically closed field and $H$ be a
finite-dimensional hereditary algebra. The associated cluster
category $\mathcal{C}_H$  was introduced and studied  in
\cite{BMRRT}, and also in \cite{CCS} for algebras $H$ of Dynkin type
$A_n$. The cluster category $\mathcal{C}_H$ is the orbit category
$D^b(H)/\tau^{-1}S$, where $S$ denotes the suspension functor and
$\tau$ is the Auslander-Reiten translation in the bounded derived
category $D^b(H)$. This is a certain 2-Calabi-Yau triangulated
category  which was invented in order to model some ingredients in
the definition of cluster algebras introduced and studied by
Fomin-Zelevinsky and Berenstein-Fomin-Zelevinsky in a series of
articles \cite{FZ1, FZ2, BFZ, FZ3}. For this purpose, a tilting
theory was developed in the cluster category. This further led to
the theory of cluster-tilted algebras initiated in \cite{BMR}.

For a positive integer $d>1$,  a certain $d+1$-Calabi-Yau category,
the $d$-cluster category $\cc_d=D^b(H)/\tau^{-1}S^d$ was considered
by Keller in \cite{Kel}. This category is showed in \cite{Tho} that
it encodes the combinatorics of the $d$-clusters of Fomin and
Reading \cite{FR} in a fashion similar to the way the cluster
category encodes the combinatorics of the clusters of Fomin and
Zelevinsky. For this reason, as a generalization of cluster
categories, the $d$-cluster category and their (cluster-)tilting
objects have been studied  in \cite{BT, IY, KR1, Tho, W, Zhu, ZZ}
and so on.

It is an interesting problem to know the algebras derived equivalent
to the cluster-tilted algebras. The study of their tilting modules
is a step in this direction. In \cite{S}, a tilting module over a
cluster-tilted algebra has been proved to lift to  a cluster-tilting
object in the  cluster category. And in \cite{FL}, the authors prove
that this result holds generally in  the 2-Calabi-Yau triangulated
category, that is, a tilting module over the endomorphism algebra of
a cluster-tilting object in a 2-Calabi-Yau triangulated category
lifts to a cluster-tilting object in this 2-Calabi-Yau triangulated
category. The aim of current note is to get similar result of
identifying tilting modules over cluster-tilted algebras
corresponding to the higher cluster category by using the
$d+1$-Calabi-Yau property. Namely, we prove the following.
\begin{thma}For a positive integer $d>1$, let $\cc_d$ be a $d$-cluster category and $T$ a $d$-cluster-tilting object, and let $\Gamma$ be
the endomorphism algebra of $T$. Then a tilting $\Gamma$-module $L$
 lifts to a $d$-cluster-tilting object in
$\cc_d$.
\end{thma}

We point out here that the methods used in this note are different
from the ones used in the case of  2-Calabi-Yau category (where
$d=1$). Some  descriptions on the relation between cluster and
classical tilting are given in \cite{HJ}.

\vspace{0.2cm} \noindent{\bf Acknowledgments.} The author would like
to thank Changjian Fu for helpful discussions. He is grateful to
Idun Reiten and Bin Zhu for valuable suggestions and comments.
Thanks also to David Smith for his interest.

\section{Preliminaries}

We recall some notations and known results.
\subsection{Tilting modules} Let $k$ be an algebraically closed field and $A$ be a finite-dimensional algebra.
Let $\mod A$ be the category of finite-dimensional right
$A$-modules. For an $A$-module $T$, let $\add T$ denote the full
subcategory of $\mod A$ with objects all direct summands of direct
sums of copies of $T$. Then $T$ is called a tilting module in $\mod
A$ if
\begin{enumerate}
\item[-]$\pd_AT\leq 1$,
\item[-]$\Ext^1_A(T,T)=0$,
\item[-]there is an exact sequence $0\to A\to T^0\to T^1\to 0$, with
$T^0, T^1$ in $\add T$.
\end{enumerate}
This is the original definition of tilting modules from \cite{HR},
and it was proved in \cite{B} that the third axiom can be replaced
by the following:
\begin{enumerate}
\item[-]the number of indecomposable direct summands of $T$ (up to
isomorphism) is the same as the number of  simple $A$-modules.
\end{enumerate}

\subsection{ $d$-Calabi-Yau categories and higher cluster categories}

Let $k$ be an algebraically closed field and $\cc$ be a
Krull-Schmidt triangulated $k$-linear category with split
idempotents and suspension functor $S$. We suppose that all
$\Hom$-spaces of $\cc$ are finite-dimensional and that $\cc$ admits
a Serre functor $\Sigma$, $cf.$\cite{RV}. Let $i\geqslant 1$ be an
integer. When we say that $\cc$ is Calabi-Yau of CY-dimension $i$
(or simply $i$-Calabi-Yau), we mean that there is an isomorphism of
triangle functors
\[S^i\stackrel{\sim}\rightarrow \Sigma.\]For $X,Y\in \cc$ and
$n\in \mathbb{Z}$, we put as usual $$\Ext_\mathcal
{C}^n(X,Y)=\Hom_\mathcal {C}(X,S^nY).$$

Let $H$ be a hereditary algebra and the number of  simple
$H$-modules be $n$ . Let $\cd=D^b(H)$ be the bounded derived
category of $H$ with suspension functor $S$ and the Auslander-Reiten
translate $\tau$. For a positive integer $d>1$, the higher cluster
category, $d$-cluster category is the orbit category
$\cc_d=\cd/\tau^{-1}S^d$. It is shown in \cite{Kel} that $\cc_d$ is
a triangulated category and the canonical functor $\cd\to \cc$ is a
triangle functor. We denote therefore by $S$ the suspension in
$\cc_d$. The $d$-cluster category is also Krull-Schmidt and is
Calabi-Yau of CY-dimension $d+1$.  That is, for any $X, Y$ in $
\cc_d,$
\[\Hom(X, Y)\simeq D\Hom(Y, S^{d+1}X),
\]or equivalently
\[\Ext^1(X, Y)\simeq D\Ext^d(Y, X).
\]

We recall the notation of $d$-cluster-tilting object  from
\cite{KR1, Tho, Zhu, IY}. This notation shares the same meaning as
"maximal $d$-orthogonal subcategory" in the sense of Iyama \cite{I}.
Let $\cc$ be a $d+1$-Calabi-Yau category. An object $X$ in $\cc$ is
called rigid if
\[\Ext^i_{\cc}(X, X)=0, \text{\quad for all}\ 1\leqslant i
\leqslant d.
\]A
rigid object $T$ is called $d$-cluster-tilting if it satisfies the
property: if $X\in \cc_d$ satisfies $\Ext^i_{\cc}(X, T)=0$ for all
$1\leqslant i\leqslant d,$ then $X\in \ct=\add T$.

Let $\cc$ be a $d+1$-Calabi-Yau category with a $d$-cluster-tilting
object $T$. Let $\Gamma$ be the endomorphism algebra of $T$. For
classes $\cu, \cv$ of objects, we denote by $\cu*\cv$ the full
subcategory of all objects $X$ of $\cc$ appearing in a triangle
\[ U\to X\to V\to SU.
\]
Let $F: \cc\to \mod \Gamma$ be the functor which sends $X$ to
$\Hom_\cc(T, X)$. There is an essential result in \cite{KR1} as
following.
\begin{theorem}For each module $M\in \mod \Gamma$,
there is a triangle
\[T_0\to T_1\to X\to ST_0
\]such that $FX$ is isomorphic to $M$. The functor $F$ induces an
equivalence
\[\ct*S\ct/(S\ct)\longrightarrow \mod \Gamma.
\]
\end{theorem}

Thus when we say that a $\Gamma$-module $L$ lifts to the higher
cluster category, we  mean its preimage under the equivalence.

We need the following theorem which is shown in \cite{Zhu}.
\begin{theorem}\label{Zhu}
Let $X$ be a rigid object in the $d$-cluster category $\cc_d$, then
$X$ is a $d$-cluster-tilting object if and only if $X$ has $n$
indecomposable summands, up to isomorphism.
\end{theorem}

\section{Proof of the main result}

First we prove the following  crucial proposition.
\begin{proposition} Let $\cc$ be a $d+1$-Calabi-Yau category with a
$d$-cluster-tilting object $T$ and let $\Gamma$ be the endomorphism
algebra of $T$. Let $M, N$ be two objects in $\ct*S\ct$ and $FM,
FN\in \mod \Gamma$ be their images under the functor $F$. If $FM$
and $FN$ are of projective dimension at most one and satisfy
$\Ext^1_\Gamma(FM, FN)=0$ and $\Ext_\Gamma^1(FN, FM)=0,$ then
$\Ext_\cc^i(M, N)=0$ and $\Ext_\cc^i(N, M)=0$ for all $0< i< d+1$.
\end{proposition}

\begin{proof}

We only need to show that the result holds for $M, N$
indecomposable. Since $FM$ is of projective dimension at most $1$,
we have the following exact sequence
\[0\to P_1^M\to P_0^M\to FM \to 0.
\]As the assumption,  $\Ext^1_\Gamma(FM, FN)=0$. By the definition of $\Ext^1,$ we
have the following commutative diagram in $\mod \Gamma$.
\[\xymatrix{P_1^M \ar[d]_{\forall f}\ar[r]& P_0^M\ar[r]\ar@{.>}[dl]^{\exists g}_{\curvearrowright} & FM\ar[r]& 0.\\
FN}
\]Using the equivalence $\ct*S\ct/(ST)\xrightarrow{\sim}\mod \Gamma$
and because $\Hom_\cc(\ct, S\ct)=0$, we have the following in $\cc$.
\[\xymatrix{T_1^M \ar[d]_{\forall \tilde{f}}\ar[r]& T_0^M\ar[r]\ar@{.>}[dl]^{\exists \tilde{g}}_{\curvearrowright} & M\ar[r]& ST_1^M.\\
N}
\]That is, for any $f: T_1^M\to N$, there exists $g: T_0^M\to N $
such that $f$ factors through $g$.

First, we claim that
\[\Hom_\cc(M, SN)=0.
\]In fact, consider the following two triangles
\[T_1^M\to T^M_0\xrightarrow{p^M_0} M\to ST^M_1
\]and
\[ST^N_1\to ST^N_0\to SN\xrightarrow{\omega} S^2T^N_1.
\]Let $\alpha$ be any morphism from $M$ to $SN$. Since $T$ is a
$d$-cluster-tilting object in $\cc$, the composation \[ \omega\cdot
\alpha\cdot p^M_0\in\Hom_\cc(T^M_0, S^2T^N_1)=0.\] Therefore there
exists a morphism from $T^M_0$ to $ST^N_0$ which makes the following
 diagram of triangles commutative.
\[\xymatrix{T^M_1\ar[r]\ar@{.>}[d]& T^M_0\ar[r]^{p^M_0}\ar@{.>}[d]& M\ar[d]^\alpha
\ar[r] & ST^M_1\\
ST^N_1\ar[r] & ST^N_0\ar[r]& SN\ar[r] & S^2T^N_1}
\]Thus we get
\[\alpha\cdot p^M_0=0
\]for the reason that $\Hom_\cc(T^M_0, ST^N_0)=0$.
So there exist $\beta: ST^M_1\to SN$ such that $\alpha$ factors
through $\beta$. As described above, we get  $\gamma: ST^M_0\to SN$
such that \[ \beta= \gamma\cdot h.\] That is, we have the following
commutative diagram
\[\xymatrix{T^M_1\ar[r]\ar@{.>}[d]& T^M_0\ar[r]^{p^M_0}\ar@{.>}[d]& M\ar[d]^\alpha
\ar[r]^u & ST^M_1\ar@{.>}[dl]_\beta\ar[r]^h & ST^M_0\ar@{.>}[dll]_\gamma\\
ST^N_1\ar[r] & ST^N_0\ar[r]& SN\ar[r] & S^2T^N_1}
\]where \[\alpha=\beta\cdot u=\gamma\cdot h\cdot u=0.\] This implies $\Hom_\cc(M, SN)=0.$
Dually one can prove that $\Hom_{\cc}(N, SM)=0.$

Now let $1<i<d$. As before, suppose that $\alpha: M\to S^iN$. Here
we consider  the following two triangles
\[T_1^M\to T^M_0\xrightarrow{p_0^M} M\to ST^M_1
\]and
\[S^iT^N_1\to S^iT^N_0\to S^iN\xrightarrow{\omega} S^{i+1}T^N_1.
\] Again because  $T$ is a $d$-cluster-tilting object and
$i+1\leqslant d$, we have  \[ \omega\cdot \alpha\cdot p^M_0\in
\Hom_\cc (T^M_0, S^{i+1}T^N_1)=0.\] Thus there exists morphism from
$T^M_0$ to $S^iT^N_0$, which is zero for the same reason that $T$ is
a $d$-cluster-tilting object, makes the diagram commutative.
\[\xymatrix{T^M_1 \ar[r]\ar@{.>}[d] & T^M_0 \ar[r]^{p^M_0}\ar@{.>}[d] &
M\ar[r] \ar[d]^\alpha & ST^M_1\ar@{.>}[dl]_\beta \\
S^iT^N_1 \ar[r]& S^iT^N_0\ar[r] & S^iN\ar[r]^\omega& S^{i+1}T^N_1
}\] That is $ \alpha\cdot p^M_0=0$, so there exists $\beta:
ST^M_1\to S^iN$ such that $\alpha$ factors through $\beta$. Note
that \[ \omega\cdot \beta\in \Hom_\cc(ST^M_1, S^{i+1}T^N_1),\] which
is zero since $T$ is a $d$-cluster-tilting object and $i>1$. So
there exists $\gamma: ST^M_1\to S^iT^N_0$ such that $\beta$ factors
through $\gamma$. But \[\Hom_\cc(ST^M_1, S^iT^N_0)=0\] for the
reason that $T$ is a $d$-cluster tilting object and $i>1$. This is
to say that $\gamma$ is zero, further more $\beta$ is zero and so is
$\alpha$.
\[\xymatrix{T^M_1 \ar[r]\ar@{.>}[d] & T^M_0 \ar[r]^{p^M_0}\ar@{.>}[d] &
M\ar[r] \ar[d]_\alpha & ST^M_1\ar@{.>}[dl]^\beta \ar@{.>}[dll]^-\gamma\\
S^iT^N_1 \ar[r]& S^iT^N_0\ar[r] & S^iN\ar[r]^\omega& S^{i+1}T^N_1
}\] Thus $\Hom_\cc(M, S^iN)=0$.

When $i=d$, thanks to the $d+1$-CY property, we have \[\Hom_\cc(M,
S^dN)=D\Hom_\cc(N, SM)=0.\]

Dually one can prove $\Ext_\cc^i(N, M)=0$ for all $0< i< d+1$.
\end{proof}

Now our main result is an easy corollary.
\begin{thma}For a positive integer $d>1$, let $\cc_d$ be a $d$-cluster category and $T$ a $d$-cluster-tilting object, and let $\Gamma$ be
the endomorphism algebra of $T$. Then a tilting $\Gamma$-module $L$
 lifts to a $d$-cluster-tilting object in
$\cc_d$.
\end{thma}

\begin{proof}The tilting $\Gamma$-module $L$ lifts to a rigid object
in $\cc_d$ by the proposition above. Note that the number of
indecomposable direct summands of $L$ (up to isomorphism) is the
same as the number of  simple $\Gamma$-modules or equivalently the
number of indecomposable summands of $T$, which is $n$ by Theorem
\ref{Zhu}. Thus $L$ lifts to a $d$-cluster-tilting object by Theorem
\ref{Zhu} again.
\end{proof}

\def\cprime{$'$}
\providecommand{\bysame}{\leavevmode\hbox
to3em{\hrulefill}\thinspace}
\providecommand{\MR}{\relax\ifhmode\unskip\space\fi MR }
\providecommand{\MRhref}[2]{%
  \href{http://www.ams.org/mathscinet-getitem?mr=#1}{#2}
} \providecommand{\href}[2]{#2}

\end{document}